\documentclass[12pt,a4paper]{amsart}

\usepackage{amssymb,amscd,latexsym,amsxtra,dsfont,enumerate}

\usepackage[markup=underlined]{changes}
\definechangesauthor[color=red]{Curt}
\definechangesauthor[color=blue]{Emanuel}
\usepackage{amsfonts,bbold}
\usepackage[mathscr]{eucal}
\usepackage{calrsfs}
\usepackage{fge}
\usepackage{pst-all}
\usepackage{mathtools}
\usepackage{scalerel}
\usepackage[autosize]{dot2texi}

\definecolor{britishracinggreen}{rgb}{0.0, 0.26, 0.15}
\definecolor{amaranth}{rgb}{0.9, 0.17, 0.31}

\newtheorem{thm}{Theorem}

\newtheorem{lem}[thm]{Lemma}
\newtheorem{prop}[thm]{Proposition}

\theoremstyle{definition}

\newcommand{\sss}[1]{{\scriptscriptstyle #1}}

\newcommand{\am}{\mathcal A}

\newcommand{\N}{\mathds N}
\newcommand{\R}{\mathds R}

\newcommand\restr[2]{{
		\left.\kern-\nulldelimiterspace 
		#1 
		\vphantom{\big|} 
		\right|_{#2} 
}}

\newcommand{\mynorm}[1]{ \left\| #1 \right\| }

\begin{document}
	
	\title{Every  symmetric Kubo-Ando connection has the order-determining property}
	
	\author{Emmanuel Chetcuti}
	\address{
		Emmanuel Chetcuti,
		Department of Mathematics\\
		Faculty of Science\\
		University of Malta\\
		Msida MSD 2080  Malta} \email {emanuel.chetcuti@um.edu.mt}
	
	\author{Curt Healey}
	\address{Curt Healey\\
		Department of Mathematics\\
		Faculty of Science\\
		University of Malta\\
		Msida MSD 2080  Malta}
	\email{curt.c.healey.13@um.edu.mt}

	\date{\today}
	\begin{abstract}
		In \cite{molnar} L.~Moln{\'a}r studied the question of whether the L\"owner partial order on the positive cone of an operator algebra is determined by the norm of any arbitrary Kubo-Ando mean. The question was affirmatively answered for certain classes of Kubo-Ando means, yet the general case was left as an open problem. We here give a complete answer to this question, by showing that the norm of every symmetric Kubo-Ando mean  is order-determining, i.e.  if $A,B\in  \mathcal B(H)^{++}$ satisfy $\Vert A\sigma X\Vert \le \Vert B\sigma X\Vert$ for every $X\in  \mathcal{A}^{{++}}$,  where $\mathcal A$ is the C*-subalgebra generated by $B-A$ and $I$, then $A\le B$.
	\end{abstract}
	\maketitle

	\section{Introduction}
	
	Recently, in \cite{molnar} the author studied the question of when the norm of a given mean, on the positive cone of an operator algebra, determines the L\"owner order.  As explained clearly in the introduction by the author, this problem is of relevance to the study of maps between positive cones of operator algebras that preserve a given norm of a given operator mean.  Such a study has received considerable attention, as can be seen for example in \cite{MR2517859, MR2836790,  MR3997630, MR4131516}.  The motivation of such investigations comes, first, from the study of norm additive maps or spectrally multiplicative maps, and secondly, from the study of the structure of certain quantum mechanical symmetry transformations relating to divergences.
	
Let us recall that a binary operation $\sigma$ on $\mathcal B(H)^{++}$ is called a \emph{Kubo-Ando connection}
if it satisfies the following properties:
	\begin{enumerate}[{\rm(i)}]
		\item If $A \leq C$ and $B \leq D$, then $ A \sigma B \leq C \sigma D $.
		\item $ C (A \sigma B) C \leq (CAC) \sigma (CBC) $.
		\item  If $ A_n  \downarrow A$ and $ B_n \downarrow B $, then $ A_n \sigma B_n \downarrow A \sigma B$.
	\end{enumerate}
	A \emph{Kubo-Ando mean} is a Kubo-Ando connection with the normalization condition $I\sigma I=I$.  The most fundamental connections are:
	\begin{itemize}
		\item  the \emph{sum} $(A,B)\mapsto A+B$,
		\item the \emph{parallel sum} $(A,B)\mapsto A:B=\left(A^{-1}+B^{-1}\right)^{-1}$,
		\item the \emph{geometric mean}
		\[(A,B)\mapsto A\sharp B=A^{\sss{\frac{1}{2}}}\,\left(A^{\sss{-\frac{1}{2}}}BA^{\sss{-\frac{1}{2}}}\right)^{\sss{\frac{1}{2}}}A^{\sss{\frac{1}{2}}}.\]
	\end{itemize}

A function $f:(0,\infty)\to(0,\infty)$ is said to be \emph{operator monotone} if $\sum_{i=1}^nf(a_i)P_i\le\sum_{j=1}^mf(b_j)Q_j$ whenever $\sum_{i=1}^na_iP_i\le\sum_{j=1}^mb_jQ_j$, where $a_i,b_j>0$, and the projections $P_i$, $Q_j$ satisfy $\sum_{i=1}^nP_i=\sum_{j=1}^mQ_j=I$.  Such a function is automatically continuous, monotonic increasing and concave.  For an operator-monotone function $f$, one has $f(A)\le f(B)$ whenever $A,B\in B(H)^{++}$ and satisfy $A\le B$.  It is easy to see that the class of operator monotone functions is closed under addition and multiplication by positive real numbers.  The transpose $f^\circ$ of the operator monotone function $f$, defined by $f^\circ(x):=xf(x^{-1})$, is again operator monotone.

Let $\sigma$ be a Kubo-Ando connection on $B(H)^{++}$.  In the proof of \cite[Lemma 3.2]{kubo} it is shown that the function $f$ defined on $(0,\infty)$ by $x\mapsto I\sigma x$ is scalar-valued, operator-monotone, and satisfies $f(B)=I\sigma B$ for every $B\in B(H)^{++}$.  This gives
\[A\sigma B=A^{\sss{\frac{1}{2}}}\,f\left(A^{\sss{-\frac{1}{2}}}BA^{\sss{-\frac{1}{2}}}\right)A^{\sss{\frac{1}{2}}}\] for every $A,B\in B(H)^{++}$.  The function $f$ is called the \emph{representing function} of $\sigma$.

We recall that operator monotone functions correspond to positive finite Borel measures on $[0,\infty]$ by L\"owner's Theorem (see \cite{MR0486556}): To every  operator monotone function $f$ corresponds a unique positive and finite Borel measure $m$ on $[0,\infty]$ such that
	\begin{equation}\label{e4}
		f(x)\,:=\,\int_{[0,\infty]}\frac{x(1+t)}{x+t}\,{\rm d}m(t)\,=\,m(\{0\})\,+\,x\,m(\{\infty\})\,+\,\int_{(0,\infty)}\frac{1+t}{t}(t:x)\,{\rm d}m(t)\quad(x>0).
	\end{equation}
It is easy to see that $f(0+)=m(\{0\})$, $f^\circ(0+)=m(\{\infty\})$.

Let $f:(0,\infty)\to(0,\infty)$ be an operator monotone function, and let $m$ be the positive and finite Borel measure on $[0,\infty]$ associated to $f$ via L\"owner's Theorem by (\ref{e4}).  The binary operation $\sigma_f$ defined  on $B(H)^{++}$ by
\[A\sigma_f B\,= \,f(0+)\, A+f^\circ(0+)\,B\,+\,\int_{(0,\infty)}\frac{1+t}{t}(tA:B)\,{\rm d}m(t)\]
satisfies conditions {\rm(i)} and {\rm(ii)} of the definition of a Kubo-Ando connection.  Moreover, 
\[I\sigma_f A\,= \,f(0+)\,I+f^\circ(0+)\,A\,+\,\int_{(0,\infty)}\frac{1+t}{t}(tI:A)\,{\rm d}m(t)\,=\,\int_{[0,\infty]}A(1+t)(t{I}+A)^{-1}\,{\rm d}m(t)\,=\,f(A)\]
for every $A\in B(H)^{++}$ and therefore
\[A\sigma_f B=A^{\sss{\frac{1}{2}}}\,f\left(A^{\sss{-\frac{1}{2}}}BA^{\sss{-\frac{1}{2}}}\right)A^{\sss{\frac{1}{2}}}\]
 for every $A,B\in B(H)^{++}$.  Using the fact that a continuous real-valued function is SOT continuous on bounded sets of self-adjoint operators (see \cite[Prop. 5.3.2, p. 327]{KadRinVol1}), it follows that if $A_n\downarrow A$ and $B_n\downarrow B$ in $B(H)^{++}$, then $A_n\sigma_f B_n\downarrow A\sigma_f B$.  This shows that $\sigma_f$ is a Kubo-Ando connection on $B(H)^{++}$.

We further recall that if $\sigma$ is a Kubo-Ando connection with representing function $f$, then the representing function of the `reversed' Kubo-Ando connection $(A,B)\mapsto B\sigma A$ is the transpose $f^\circ$.  The Kubo-Ando connection is said to be symmetric if it coincides with its reverse, i.e. a Kub-Ando connection is symmetric if and only if the representing function $f$ satisfies $f=f^\circ$ as shown in \cite[Corollary 4.2]{kubo}. The Kubo-Ando means are precisely the Kubo-Ando connections whose representing function satisfy the normalizing condition $f(1)=1$.
	
The most fundamental Kubo-Ando means are the power means which correspond to the operator monotone functions
	\[f_{\sss{p}}(t):=\begin{cases}
		\left(\dfrac{1+t^p}{2}\right)^{\frac{1}{p}}, & \mbox{if } -1\le p\le 1,\ p\neq 0 \\
		\sqrt{t}, & \mbox{if } p=0.
	\end{cases}\]
	The principal cases $f_{\sss{0}}(t)=\sqrt{t}$, $f_{\sss{-1}}=\frac{2t}{1+t}$ and $f_{\sss{1}}(t)=\frac{1+t}{2}$ correspond, respectively, to the geometric mean $(A,B)\mapsto A\sharp B$, the harmonic mean $(A,B)\mapsto A!B=2(A:B)$, and the arithmetic mean $(A,B)\mapsto A\nabla B=(A+B)/2$.
	
It is easy to verify that the measure $m$ associated to the arithmetic mean is $(\delta_{\sss{ 0}}+\delta_{\sss{ \infty}})/2$ and that associated to the harmonic mean is $\delta_{\sss{ 1}}$, where $\delta_{\sss{ x}}$  denotes the Dirac measure on the point $x\in[0,\infty]$.

\emph{A remark on the domain of definition of a Kubo-Ando connection:}  We have opted for having $B(H)^{++}$ the defining domain of a Kubo-Ando connection (as opposed to $B(H)^+$) in order to obtain  fully consistent interchangeable relations in the diagram below.
\begin{figure}
\begin{center}
\begin{pspicture}(-5,-4)(6,6)
\pnode(-5,-2){a}
\pnode(5,-2){b}
\pnode(0,5){c}
\uput[l](a){$f$}
\uput[r](b){$m$}
\uput[d](c){$\sigma$}
\uput[d](a){${\scriptstyle \text{Operator monotone function } f}$}
\uput[d](b){${\scriptstyle \text{Positive and finite Borel measure } m \text{ on }[0,\infty]}$}
\uput[u](c){${\scriptstyle \text{Kubo-Ando connection } \sigma \text{ on }B(H)^{++}}$}
\psline[linewidth=0.5pt]{<->}(-4.8,-2)(4.8,-2)
\psline[linewidth=0.5pt]{<->}(-4.8,-1.8)(-0.1,4.6)
\psline[linewidth=0.5pt]{<->}(4.8,-1.8)(0.1,4.6)
\uput[l](-2.5,1.5){${\scriptstyle A\sigma B=A^{\sss{\frac{1}{2}}}\,f\left(A^{\sss{-\frac{1}{2}}}BA^{\sss{-\frac{1}{2}}}\right)A^{\sss{\frac{1}{2}}}}$}
\uput[r](2.5,1.5){${\scriptstyle A\sigma B=m(\{0\})\, A+m(\{\infty\})\,B+\int_{(0,\infty)}\frac{1+t}{t}(tA:B)\,{\rm d}m(t) }$}
\uput[d](0,-2.5){${\scriptstyle f(x)=\int_{[0,\infty]}\frac{x(1+t)}{x+t}\,{\rm d}m(t)}$}
\end{pspicture}
\end{center}
\end{figure}
 It must be said, however, that this offers no handicap because any Kubo-Ando connection $\sigma$ on $B(H)^{++}$ can be extended to a binary relation $\hat\sigma$ on $B(H)^+$ by setting $A\hat\sigma B$ equal to
\[\inf\{(A+n^{-1})\sigma (B+n^{-1}):n\in\N\}=\inf\{X\sigma Y:X,Y\in B(H)^{++}, X\ge A, Y\ge B\},\]
and it is not hard to show that the extension $\hat\sigma$ satisfies {\rm(i)}-{\rm(iii)} of the definition of a Kubo-Ando connection.  Note that in this case the equality
\[A\hat\sigma B=A^{\sss{\frac{1}{2}}}\,f\left(A^{\sss{-\frac{1}{2}}}BA^{\sss{-\frac{1}{2}}}\right)A^{\sss{\frac{1}{2}}},\]
where $f$ is the representing function associated to $\sigma$, holds only on $B(H)^{++}\times B(H)^{+}$.  We further remark that the continuity properties of the function calculus (see \cite[Prop. 5.3.2, p. 327]{KadRinVol1}) imply that:
\begin{enumerate}[{\rm(i)}]
\item The map
\begin{equation}\label{e8}\hat\sigma\,:\,\mathcal B(H)^{++}\times \mathcal B(H)^{+}\,\to\, \mathcal B(H)^{+}\,:\,(A,B)\,\mapsto\, A^{\sss{\frac{1}{2}}}\,f\left(A^{\sss{-\frac{1}{2}}}BA^{\sss{-\frac{1}{2}}}\right)A^{\sss{\frac{1}{2}}}\end{equation}
is continuous when the domain is equipped with the product of the relative topologies induced by the norm, and the range with the norm topology, and
\item for every $\epsilon, R>0$, the restriction of $\hat\sigma$ to the rectangle
\[ [\epsilon I\,,\,RI]\times[0\,,\,RI]\]
is continuous when the domain is equipped with the product of the relative topologies induced by SOT, and the range with SOT.
\end{enumerate}
In the sequel we will not distinguish between $\sigma$ and $\hat\sigma$ any further.

\section{Preliminary Considerations}
	
In this section, we collate a list of lemmas and propositions which will prove to be helpful in proving the main result.
	
\begin{lem}\label{l1}
Let $f:(0, \infty)  \rightarrow (0, \infty) $ be an operator monotone  function and let $m$ denote the positive and finite Borel measure associated to $f$ via (\ref{e4}).
\begin{enumerate}[{\rm(i)}]
\item  If $f$ is symmetric,  $\int_{[0,\infty]}t\,{\rm d}m(t)=\int_{[0,\infty]}t^{-1}\,{\rm d}m(t)$.

\item For every Borel subset $\Delta$ of $[0,\infty]$, the function $f_\Delta$ defined on $(0,\infty)$  by
\[f_\Delta:x\mapsto\int_{\Delta}\frac{x(1+t)}{x+t}\,{\rm d}m(t)\]
is operator monotone.  In particular, the function $h$ defined by
\[h(x):=\int_{(0,\infty)}\frac{x(1+t)}{x+t}\,{\rm d}m(t)=f(x)-f(0+)-f^\circ(0+) x\quad(x>0)\]
is operator monotone.  If $f$ is symmetric, then so is $h$.
\end{enumerate}
	\end{lem}
\begin{proof}
{\rm(i)}~By the Monotone Convergence Theorem one has that
\[f(x)=\int_{[0,\infty]} \frac{1+t}{1+tx^{-1}}\,{\rm d}m(t)\,\uparrow\, \int_{[0,\infty]}1+t\,{\rm d}m(t)\qquad \text{as $x\uparrow\infty$,}\]
 and since $f^\circ(x)=xf(1/x)$ one also gets
\[f^\circ(x)=\int_{[0,\infty]}\frac{1+t}{t+x^{-1}}\,{\rm d}m(t)\, \uparrow\, \int_{[0,\infty]}1+t^{-1}\,{\rm d}m(t)\text{as $x\uparrow\infty$.}\]
So, if $f$ is symmetric, $\int_{[0,\infty]}t\,{\rm d}m(t)=\int_{[0,\infty]}t^{-1}\,{\rm d}m(t)$.

{\rm(ii)}~Since
\[\infty>f(x)=\int_{[0,\infty]}\frac{x(1+t)}{x+t}\,{\rm d}m(t)\ge \int_{\Delta}\frac{x(1+t)}{x+t}\,{\rm d}m(t),\]
the Lebesgue Dominated Convergence Theorem can be applied to deduce that the function $f_\Delta$ is continuous.

The function
\[\mathcal B(H)^{++}\ni A\mapsto f_\Delta(A)=\int_\Delta (1+t)(1+tA^{-1})^{-1}\,{\rm d}m(t)\]
maps $\mathcal B(H)^{++}$ into $\mathcal B(H)^{+}$ and satisfies $f_\Delta(A)\le f_\Delta(B)$ whenever $A\le B$ in $\mathcal B(H)^{++}$.  Given $A\le B$ in $B(H)^+$, the operators $A_n:=A+n^{-1}$, $B_n:=B+n^{-1}$ ($n\in\N$) belong to $B(H)^{++}$ and satisfy $A_n\le B_n$ for every $n$.  Therefore, $f_\Delta(A_n)\le f_\Delta(B_n)$ for every $n$.  The continuity of $f_\Delta$ (see \cite[Prop. 5.3.2, p. 327]{KadRinVol1}) yields the required inequality $f_\Delta(A)\le f_\Delta(B)$.

Setting $\Delta=(0,\infty)$ one obtains that
\[h(x):=\int_{(0,\infty)}\frac{x(1+t)}{x+t}\,{\rm d}m(t)\quad(x>0)\]
is operator monotone.  Since $m(\{0\})=f(0+)$ and $m(\{\infty\})=f^\circ(0+)$,
\[h(x)=f(x)-f(0+)-f^\circ(0+)x\]
follows.  It is easy to verify that if $f$ is symmetric, so is $h$.
\end{proof}

By function calculus, for an arbitrary self-adjoint $X\in \mathcal B(H)$,  there exists a $\ast$-isomorphism from the Banach algebra  $L^{\infty}(\Delta,B(\Delta),\mu)$, where $\Delta$ is the spectrum of $X$, $B(\Delta)$ is the Borel $\sigma$-algebra in $\Delta$, and $\mu$ is the composition of the trace function with the spectral measure associated to $X$, into  the von Neumann subalgebra of $\mathcal B(H)$ generated by $X$ and $I$.  This $\ast$-isomorphism preserves suprema/infima of monotone sequences and maps the continuous functions on $\Delta$ isometrically onto  the C*-subalgebra of $\mathcal B(H)$ generated by $X$ and $I$.  Denote by $\Gamma(X)$ the range of this $\ast$-isomorphism.

\begin{lem}\label{lemma_using_spectral_projections}
Let $\sigma$ be a Kubo-Ando connection.  For $A,B\in \mathcal B(H)^{++}$ let $\mathcal A$ denote the  C*-subalgebra of $\mathcal B(H)$ generated by $B-A$ and $I$.   If $\mynorm{A \sigma X} \leq \mynorm{B \sigma X} $ holds for all $ X \in \mathcal{A}^{++} $, then $\mynorm{A \sigma X} \leq \mynorm{B \sigma X}$ holds for all positive $X\in \Gamma(B-A)$.
\end{lem}
\begin{proof}
Let us first recall the general fact that whenever $(X_\gamma)$ is an SOT-convergent net of positive operators, bounded from above by its SOT-limit $X$, then the net of norms $\left(\Vert X_\gamma\Vert\right)$ is convergent to $\Vert X\Vert$.

Invoking the continuity of $\sigma$ w.r.t. the norm (as mentioned in the introduction), and since $\mathcal A^{++}$ is norm-dense in $\mathcal A^+$, it can be seen that the hypothesis implies that $\mynorm{A\sigma X}\le\mynorm{B\sigma X}$ holds for every $X\in \mathcal A^+$.

Moreover, since for every positive Borel function on $\Delta$, there exists a sequence $(h_n)$ of positive continuous functions on $\Delta$ satisfying $h_n(x)\uparrow f(x)$ for every $x\in\Delta$,  for every $X\in\Gamma(B-A)$, one can find  a sequence $(X_n)$ in $\am^{+}$ with $X_n\uparrow X$.  Invoking now the monotonicity property and the SOT-continuity of $\sigma$ (as mentioned in the introduction) it can be deduced that $A\sigma X_n\uparrow A\sigma X$ and $B\sigma X_n\uparrow B\sigma X$.   So, the assertion now follows by the recall announced in the beginning of the proof.
\end{proof}

In the subsequent lemma, the main ideas can be found in \cite[Lemma 11]{chabbabi}. We formalise them and present them here for completeness sake.
	
\begin{prop}\label{p2}
For the operators $ A, B \in \mathcal B(H)^+$  the following assertions are equivalent:
\begin{enumerate}[{\rm(i)}]
\item  $A\le B$,
\item $\mynorm{PAP} \leq \mynorm{PBP} $ for every spectral projection $P$ of $B-A$,
\item $\{\lambda\ge 0:\lambda P\le PAP\}\subseteq \{\lambda\ge 0:\lambda P\le PBP\}$ for every spectral projection $P$ of $B-A$.
\end{enumerate}
\end{prop}
	\begin{proof}
The assertions {\rm(i)}$\Rightarrow ${\rm(ii)} and {\rm(i)}$\Rightarrow ${\rm(iii)} are trivial.

Suppose that $A\nleq B$. Then, there exists $\varepsilon > 0 $ such that the spectrum of $B-A$ has a nontrivial intersection with $(-\infty, -\epsilon)$. Let $\Delta:=(-\infty,-\epsilon)$ and let $ P_\epsilon $ be the  (non-zero) spectral projection  of $B-A$ associated to the indicator function $ \chi_{{}_\sss{\Delta}}$.  Clearly,  $ t\,\chi_{{}_\sss{\Delta}}(t)\le -\epsilon\, \chi_{{}_\sss{\Delta}}(t)$ for every $t\in\R$, and therefore $P_\epsilon BP_\epsilon-P_\epsilon AP_\epsilon=P_\epsilon (B-A)P_\epsilon \leq -\epsilon P_\epsilon$.

Rearranging the terms, we get $P_\epsilon BP_\epsilon \le P_\epsilon (A-\epsilon I) P_\epsilon$ and therefore $\mynorm{P_\epsilon BP_\epsilon} \leq \mynorm{P_\epsilon(A-\epsilon I)P_\epsilon}=\mynorm{P_\epsilon AP_\epsilon} - \epsilon$.  This shows that {\rm(ii)} implies {\rm(i)}.

Let us prove that {\rm(iii)} implies {\rm(ii)}.  First observe that for every $A\in \mathcal B(H)^+$ and projection $P$, the supremum of $\{\lambda\ge 0:\lambda P\le PAP\}$ is indeed a maximum and is at most equal to $\Vert A\Vert$.  Let $\lambda_0:=\max\{\lambda\ge 0:\lambda P_\epsilon\le P_\epsilon A P_\epsilon\}$.  Then {\rm(iii)} implies that $\lambda_0P_\epsilon\le P_\epsilon B P_\epsilon\le P_\epsilon A P_\epsilon-\epsilon P_\epsilon$, which in turn shows that $(\lambda_0+\epsilon)P_{\epsilon}\le P_\epsilon A P_\epsilon$, contradicting the maximality condition of $\lambda_0$.  This shows that {\rm(iii)} implies {\rm(i)}.
	\end{proof}
	
In the following proposition,  the ideas in \cite[Proposition 10]{molnar} are used to generalize \cite[Equation 15]{molnar}.  This will be of pivotal importance in proving the main result of this paper.
	
\begin{prop}\label{p3}
Let $X_s\in \mathcal B(H)^+$, $s>0$ satisfy $\lim_{s\to\infty}X_s=X$ in norm, and let $P\in \mathcal B(H)$ be a projection.  Then
\[  \lim_{s\to \infty} \mynorm{X_s + sP} -s = \mynorm{PXP}.\]
\end{prop}
\begin{proof}
Let $\varepsilon>0$.  It can easily be verified that
		\[\bigl\|\bigl((\mynorm{PX_sP}+\varepsilon)^{\sss{-1/2}}P+s^{\sss{-1/2}}(I-P)\bigr)\,X_s\,\bigl((\mynorm{PX_sP}+\varepsilon)^{\sss{-1/2}}P+s^{\sss{-1/2}}(I-P)\bigr)\bigr\|\]
		converges to $(\mynorm{PXP}+\varepsilon)^{-1}\mynorm{PXP}<1$ as $s\to \infty$.  Therefore,  for sufficiently large $s$
		\[\bigl((\mynorm{PX_sP}+\varepsilon)^{\sss{-1/2}}P+s^{\sss{-1/2}}(I-P)\bigr)\,X_s\,\bigl((\mynorm{PX_sP}+\varepsilon)^{\sss{-1/2}}P+s^{\sss{-1/2}}(I-P)\bigr)\ \le\ I,\]
		or $X_s\le (\mynorm{PX_sP}+\varepsilon)P+s(I-P)$.   This implies that
		\begin{equation}\label{e2}\mynorm{X_s+sP}-s\,\le\,\mynorm{PX_s P}+\varepsilon,
		\end{equation}
		for sufficiently large $s$. On the other-hand, for every $s>0$,
		\[\mynorm{X_s+sP}\ge\mynorm{PX_sP+sP},\]
		and therefore,
		\begin{equation}\label{e3}\mynorm{X_s+sP}-s\ge\mynorm{PX_sP+sP}-s=\mynorm{PX_s P}.
		\end{equation}
		Combining (\ref{e2}) and (\ref{e3}),  for sufficiently large $s$ it holds that
		\[\mynorm{PX_s P}\,\le\,\mynorm{X_s+sP}-s\,\le\,\mynorm{PX_s P}+\varepsilon.\]
		This proves that $\lim_{s\to \infty} \mynorm{X_s + sP} -s = \mynorm{PXP}$.
	\end{proof}

\begin{prop}\label{p4}\cite[Lemma 2]{molnar}
Let $f:(0, \infty)  \rightarrow (0, \infty) $ be a nontrivial (i.e. not affine) operator monotone  function satisfying $f(0+)=0$ and let $\sigma$ denote the Kubo-Ando connection associated to $f$ via (\ref{e8}).  For  $A\in \mathcal B(H)^{++}$ and nonzero projection $P\in \mathcal B(H)$
\[\mynorm{A\sigma P}=f^\circ\left(\frac{1}{\max\{\lambda\ge0:\lambda P\le PA^{-1}P\}}\right).\]
\end{prop}
	
\section{Results}

\begin{thm}\label{main theorem}
Let $\sigma$ be a  nontrivial symmetric Kubo-Ando connection on $\mathcal B(H)^{++}$. Then for every $A,B\in \mathcal B(H)^{++}$
		\[  A \leq B \iff  \mynorm{A \sigma X} \leq \mynorm{B \sigma X}, \quad\forall X \in \mathcal{A}^{++}, \]
		where $\am$ equals the  C*-subalgebra of $\mathcal B(H)$ generated by $B-A$ and $I$.
	\end{thm}
\begin{proof}  The implication $\Rightarrow$ follows trivially by the monotonicity property of Kubo-Ando means. We shall show the converse.  By Lemma \ref{lemma_using_spectral_projections} we can suppose that $\Vert A\sigma X\Vert \le \Vert B\sigma X\Vert$ holds for every $X$ of the form $X=s P+tI$, where $P$ is a spectral projection of $B-A$ and $s,t\in\R^+$.  Let $f$ be the operator monotone function associated to $\sigma$ and let $m$ be the positive and finite Borel measure associated to $f$ by L\"owner's Theorem.  Let $\alpha=f(0+)=m(\{0\})$.  The proof will be divided in cases.\\

\emph{Case 1: $\alpha=0$}.  Since $f=f^\circ$ is strictly monotonic increasing, this case follows immediately by Propositions \ref{p4} and \ref{p2}.  This result was obtained by L. Moln{\'a}r in \cite{molnar}. \\

\emph{Case 2a: $\alpha\neq 0$ and $\int_{(0,\infty)}t\,{\rm d}m(t)<\infty$}.    Let $\gamma:=\int_{(0,\infty)}1+t\,{\rm d}m(t)$.  For every $s,t,\delta>0$,  $A\in \mathcal B(H)^{++}$ and nonzero projection $P\in \mathcal B(H)$:
		\begin{multline*}
			\int_{(0,\infty)}\frac{1+t}{t}\left(tA:sP+s\delta I\right)\,{\rm d}m(t)\ -\ A\int_{(0,\infty)}1+t\,{\rm d}m(t)\\
			=\int_{(0,\infty)}\left(A\left(\left(\frac{tA}{s}+P+\delta I\right)^{-1}(P+\delta I)-I\right)\right)(1+t)\,{\rm d}m(t).
		\end{multline*}
		Noting that $\left\Vert A\left(\left(\frac{tA}{s}+P+\delta I\right)^{-1}(P+\delta I)-I\right)\right\Vert$ is a bounded function of $s$ and $t$, and using the fact that $\int_{(0,\infty)}1+t\,{\rm d}m(t)<\infty$, it is possible to apply the Dominated Convergence Theorem to infer that
		\[\int_{(0,\infty)}\frac{1+t}{t}\left(tA:sP+s\delta I\right)\,{\rm d}m(t)\]
		converges in norm to $\gamma\,A$ as $s\to\infty$.  This implies that
		\[A\sigma(sP+s\delta I)-\beta(sP+s\delta I)\ \rightarrow\  (\alpha+\gamma) A\]
		in norm, as $s\to\infty$.  Noting that $\beta=m(\{\infty\})=m(\{0\})>0$ and applying Proposition \ref{p3},  it is deduced that
		\[\lim_{s\to\infty}\bigl (\Vert A\sigma (sP+s\delta I)-\beta s\delta I\Vert-\beta s\bigr)=(\alpha+\gamma)\Vert P AP\Vert.\]
		
Using the fact that
		\[A\sigma (sP+s\delta I)=\alpha A+\beta(sP+s\delta I)+\int_{(0,\infty)}\frac{1+t}{t}(tA:sP+s\delta I)\,{\rm d}m(t)\ge \beta s\delta I,\]
		it can be seen that $\Vert A\sigma (sP+s\delta I)-\beta s\delta I\Vert=\Vert A\sigma(sP+s\delta I)\Vert-\beta s\delta$.  This establishes that
		\begin{equation}\label{e1}\lim_{s\to\infty}\bigl(\Vert A\sigma (sP+s\delta I)\Vert-\beta s(1+\delta)\bigr)=(\alpha+\gamma)\Vert P A P\Vert .\end{equation}
		
So, if $A, B\in \mathcal B(H)^{++}$ satisfy $\mynorm{A\sigma X}\le \mynorm{B\sigma X}$ for every $X$ of the form $X=s P+tI$, where $P$ is a spectral projection of $B-A$ and $s,t \in\R^+$, it follows that $\mynorm{PAP}\le\mynorm{PBP}$  holds for every spectral projection of $B-A$ and the result follows by Proposition \ref{p2}.\\

\emph{Case 2b:  $\alpha\neq 0$ and $\int _{(0,\infty)}t\,{\rm d}m(t)=\infty$}.  Denote by $\sigma_h$ the (symmetric) Kubo-Ando connection associated to the function $h(x)=f(x)-\alpha-\alpha x$  (see  {\rm(ii)} of  Lemma \ref{l1}).  Let $m_h$ denote the positive and finite Borel measure associated to $h$.  Then, $m_h(\Delta)=m(\Delta\cap (0,\infty))$ for every Borel subset $\Delta$ of $[0,\infty]$.

We suppose that
\begin{equation}\label{e9}
\Vert \alpha A+\alpha s P+A\sigma_h(sP)\Vert\le \Vert \alpha B+\alpha s P+B\sigma_h(sP)\Vert
\end{equation}
for every spectral projection $P$ of $B-A$ and $s>0$.

Noting that $h(sP)=h(s)P$ for every $s>0$, it can then be deduced that
\[A\sigma_h (sP)\le (\Vert A\Vert I)\sigma_h (sP)=\Vert A\Vert h(\Vert A\Vert^{-1} (sP))=\Vert A\Vert h(\Vert A\Vert^{-1}s)\,P\]
i.e. $A\sigma_h P$ commutes with $P$.  Thus, (\ref{e9}) yields
\begin{align*}
\Vert \alpha s P\Vert+\mynorm{ A\sigma_h(sP)} =&\mynorm{sP+A\sigma_h (sP)}\\
\leq&\Vert \alpha A+\alpha s P+A\sigma_h(sP)\Vert\\
\le& \Vert \alpha B+\alpha s P+B\sigma_h(sP)\Vert\\
\le& \Vert\alpha B\Vert+\Vert\alpha sP\Vert+\Vert B\sigma_h(sP)\Vert,
\end{align*}
and therefore
\begin{equation}\label{e10}
 \mynorm{ A\sigma_h(sP)}-\Vert B\sigma_h(sP)\Vert\le \Vert\alpha B\Vert.
\end{equation}

 Let $c_A:= 1/\max\{\lambda\ge 0:\lambda P\le PA^{-1}P\}$ and let $c_B$ be defined similarly.  Proposition \ref{p4} gives
\[\Vert A\sigma_h(sP)\Vert =s\Vert (s^{-1}A) \sigma_h P\Vert = s\, h\left(\frac{1}{s\,\max\{\lambda\ge 0:\lambda P\le PA^{-1}P\}}\right)=s\,h\left({c_A}{s^{-1}}\right)\]
and since $m_h$ is just the restriction of the measure $m$ (associated to $f$) to $(0,\infty)$, we obtain
\[\Vert A\sigma_h(sP)\Vert=s\,h\left({c_A}{s^{-1}}\right)=\int_{[0,\infty]}\frac{s c_A (1+t)}{c_A + st}\, {\rm d}m_h(t)=\int_{(0,\infty)}\frac{s c_A (1+t)}{c_A + st}\, {\rm d}m(t).\]
Similarly, $\Vert B\sigma_h(sP)\Vert=\int_{(0,\infty)}\frac{sc_A (1+t)}{c_A + st}\, {\rm d}m(t)$ and therefore,
\begin{align*}
 \mynorm{ A\sigma_h(sP)}-\Vert B\sigma_h(sP)\Vert &= \int_{(0, \infty)} \frac{sc_A (1+t)}{c_A + st} -  \frac{sc_B (1+t)}{{c_B} + st} \, {\rm d}m(t)\\
&= (c_A -c_B) \int_{(0, \infty)}\frac{s^{2}t(1+t)}{ (c_A + st)(c_B + st)} \, {\rm d}m(t).
\end{align*}

The Monotone Convergence Theorem implies that as $s\uparrow \infty $,  the integral increases to $\int_{(0, \infty)} \frac{1}{t} + 1 \, {\rm d}m(t)$.  The relation between the two measures $m$ and $m_h$,  part {\rm(i)} of  Lemma \ref{l1} and the hypothesis, then yield that
\begin{align*}
\mynorm{ A\sigma_h(sP)}-\Vert B\sigma_h(sP)\Vert\,&\uparrow\,\int_{(0, \infty)} \frac{1}{t} +1\, {\rm d}m(t)\\
&=\,\int_{[0, \infty]} \frac{1}{t}+1 \, {\rm d}m_h(t)=\int_{[0, \infty]}  t+1 \, {\rm d}m_h(t)=\int_{(0, \infty)} {t} +1 \, {\rm d}m(t)=\infty
\end{align*}
as $s\uparrow \infty$.  Therefore, $ c_{A} - c_{B} \leq 0 $ since otherwise one would get a contradiction with (\ref{e10}).  This shows that
\[ \max\{\lambda\ge 0:\lambda P\le PA^{-1}P\}\ge  \max\{\lambda\ge 0:\lambda P\le PB^{-1}P\}\]
for every spectral projection of $B-A$, and therefore $A^{-1}\ge B^{-1}$ by Proposition \ref{p2}.

\end{proof}

\subsection*{Acknowledgement}  The authors are grateful to Professor Lajos Moln{\'a}r (Bolyai Institute and University of Szeged)  who introduced the topic and problem to them.
	
	
	\bibliographystyle{plain}

\begin{thebibliography}{1}

\bibitem{chabbabi}
Fadil Chabbabi, Mostafa Mbekhta, and Lajos Moln\'{a}r.
\newblock Characterizations of {J}ordan {$^*$}-isomorphisms of {$C^*$}-algebras
  by weighted geometric mean related operations and quantities.
\newblock {\em Linear Algebra Appl.}, 588:364--390, 2020.

\bibitem{MR0486556}
William~F. Donoghue, Jr.
\newblock {\em Monotone matrix functions and analytic continuation}.
\newblock Die Grundlehren der mathematischen Wissenschaften, Band 207.
  Springer-Verlag, New York-Heidelberg, 1974.

\bibitem{KadRinVol1}
Richard~V. Kadison and John~R. Ringrose.
\newblock {\em Fundamentals of the theory of operator algebras. {V}ol. {I}},
  volume~15 of {\em Graduate Studies in Mathematics}.
\newblock American Mathematical Society, Providence, RI, 1997.
\newblock Elementary theory, Reprint of the 1983 original.

\bibitem{kubo}
Fumio Kubo and Tsuyoshi Ando.
\newblock Means of positive linear operators.
\newblock {\em Math. Ann.}, 246(3):205--224, 1979/80.

\bibitem{MR2517859}
Lajos Moln\'{a}r.
\newblock Maps preserving the harmonic mean or the parallel sum of positive
  operators.
\newblock {\em Linear Algebra Appl.}, 430(11-12):3058--3065, 2009.

\bibitem{MR2836790}
Lajos Moln\'{a}r.
\newblock Maps preserving general means of positive operators.
\newblock {\em Electron. J. Linear Algebra}, 22:864--874, 2011.

\bibitem{MR3997630}
Lajos Moln\'{a}r.
\newblock Quantum {R}\'{e}nyi relative entropies on density spaces of
  {$C^\ast$}-algebras: their symmetries and their essential difference.
\newblock {\em J. Funct. Anal.}, 277(9):3098--3130, 2019.

\bibitem{MR4131516}
Lajos Moln\'{a}r.
\newblock Maps on positive cones in operator algebras preserving power means.
\newblock {\em Aequationes Math.}, 94(4):703--722, 2020.

\bibitem{molnar}
Lajos Moln\'{a}r.
\newblock On the order determining property of the norm of a {K}ubo-{A}ndo mean
  in operator algebras.
\newblock {\em Integral Equations Operator Theory}, 93(5):Paper No. 53, 25,
  2021.

\end{thebibliography}

\end{document}